\documentclass[12pt,oneside,reqno]{amsart}
\usepackage{amssymb}
\usepackage{amsmath}
\usepackage{amsthm}
\usepackage{longtable}
\usepackage{xspace}
\usepackage{anysize} 

\allowdisplaybreaks

\theoremstyle{plain}
\newtheorem{Theorem}{Theorem}[section]

\newtheorem{Proposition}[Theorem]{Proposition}
\newtheorem{Lemma}[Theorem]{Lemma}

\newtheorem*{Remark*}{Remark}
\newtheorem*{Remarks*}{Remarks}
\newtheorem*{Theorem*}{Theorem}
\newtheorem*{Proposition*}{Proposition}
\newtheorem*{Corollary*}{Corollary}

\newcommand{\C}{\mathbb{C}}

\newcommand{\Q}{\mathbb{Q}}

\newcommand{\Z}{\mathbb{Z}}
\newcommand{\re}{\operatorname{Re}}

\newcommand{\ordp}{\mathrm{ord}_p}
\newcommand{\prodchi}{\prod_{\substack{\chi\,\mathrm{mod}\,p^{n+1}\\ \chi\,\mathrm{odd}}}}
\newcommand{\prodchieven}{\prod_{\substack{\chi\,\mathrm{mod}\,p^{n+1}\\ \chi\,\mathrm{even}}}}
\newcommand{\Mod}[1]{\mathrm{mod}\ #1}

\begin{document}

\pagestyle{plain}

\title{On the growth of the $(S,\{2\})$-refined class number}

\author{Eugenio Finat}
\address{Universidad de Chile, Facultad de Ciencias, Casilla 653, Santiago, Chile}
\email{e\_finat@yahoo.com}
\date{\today}

\begin{abstract}
We give a new proof of the fact that the growth (with respect to $n$) of the $p$-adic 
valuation of $h_{n,2}^{-}$ is linear, where 
$h_{n,2}^{-}$ denotes the minus part of the $(S,\{2\})$-refined class number of 
the cyclotomic field $\Q(\mu_{p^{n+1}})$, as defined by Hu and Kim in 
\cite{HK}. As a consequence of our proof, 
we obtain an explicit relation between the $p$-adic 
valuation of $h_{n,2}^{-}$ and the $p$-adic 
valuation of $h_{n}^{-}$,  the minus part of the class number $h_n$ of the 
cyclotomic field $\Q(\mu_{p^{n+1}})$.
\end{abstract}


\maketitle

\section{Introduction}
\label{intro}

Let $p$ be a fixed {\bf odd} prime number and, for any $x\in\Q$,
let $\ordp\,x$ be its $p$-adic valuation.
We will use the notation $n\gg0$  for ``sufficiently large $n$''.

For integers $\ell\ge1$, let $\mu_\ell\subset\C$ be the group of $\ell$-th roots of unity and 
let $\Q(\mu_\ell)$ be the $\ell$-th cyclotomic field. 
Let $h$ be the class number of $\Q(\mu_\ell)$, this is, $h$ is the
order of the ideal class group of $\Q(\mu_\ell)$  (i.e., the group of 
nonzero fractional ideals of $\Q(\mu_\ell)$
modulo its subgroup of principal fractional ideals). Then $h$ admits a factorization 
$h=h^+h^-$, where $h^+$ is the class number 
of the maximal real subfield of 
$\Q(\mu_\ell)$, and $h^-=h/h^+$ is an integer 
called the minus part of $h$ (or the relative class number).

A classical result in the arithmetic theory of cyclotomic fields, conjectured and partially
proved by Iwasawa, states the following (see \cite[Theorem 3.2 in p.~260]{LangCF}).
\begin{Theorem} 
\label{Theohn}
Let $h_n^{-}$ be the minus part of the class number $h_n$ of the 
cyclotomic field $\Q(\mu_{p^{n+1}})$. Then there exist constants 
$\lambda\ge0$ and $c$ such that
\begin{equation}
\notag 
\ordp\,h_n^{-}=\lambda n+c\qquad(n\gg0).
\end{equation}
\end{Theorem}

Iwasawa originally proved that (see \cite[Theorem 1 in p.~94]{IwaL})
\begin{equation}
\label{iwasawa}
\ordp\,h_n^{-}=mp^n+\lambda n+c\qquad(n\gg0),
\end{equation}
hence, Theorem \ref{Theohn} says that $m=0$ in \eqref{iwasawa}, this is, 
 the growth (with respect to $n$) of $\ordp\,h_n^{-}$ is 
linear. 
The proof of Theorem \ref{Theohn} involves computing the $p$-adic valuation of the 
right-hand side in the formula (see Lang \cite[p.80]{LangCF})
\begin{equation}
\label{hnbn}
h_n^{-}=2p^{n+1}\prodchi\left(-\frac12B_{1,\chi}\right).
\end{equation}
Here, the product is taken over all odd primitive Dirichlet characters of conductor 
dividing $p^{n+1}$, and $B_{1,\chi}$ is the first generalized Bernoulli number attached
to the character $\chi$, which is a special value of the Dirichlet $L$-function 
$L(s,\chi)$, namely \cite[Theorem 1 in \S2]{IwaL},
\begin{equation}
\label{lzero}
L(0,\chi)=-B_{1,\chi}.
\end{equation}

In the recent article \cite{HK}, Hu and Kim proved analogues of \eqref{iwasawa}
and Theorem \ref{Theohn} in the context
of $(S,T)$-refined class groups. 
For sake of brevity
we will not enter into the details of this theory
and we refer the reader to \cite[\S2]{HK} and to the references mentioned therein.
Let $\chi$ be a Dirichlet character and let $L_E(s,\chi)$ be the alternating (or Eulerian) Dirichlet 
$L$-function defined by the series
\begin{equation}
\label{LEdef}
L_E(s,\chi)=2\sum_{k\ge1}(-1)^k\frac{\chi(k)}{k^s}\qquad(\re(s)>0),
\end{equation}
which can be analytically continued to the entire complex plane. 
Let $S$ be the set
of infinite places of $\Q(\mu_{p^{n+1}})$ and let $T$ 
be the set of places above the rational prime $2$ in the same field.
If $h_{n,2}^{-}$ denotes the minus part of the $(S,T)$-refined class number of 
$\Q(\mu_{p^{n+1}})$, then Hu and Kim proved that (see \cite[Proposition 3.4]{HK}) 
\begin{equation}
\label{hn2def}
h_{n,2}^{-}=(-1)^{\varphi(p^{n+1})/2}\,2^{1-\varphi(p^{n+1})}
\prodchi L_E(0,\chi),
\end{equation}
in analogy with \eqref{hnbn}. Here, $\varphi$ is Euler's totient function. 
Using some results in \cite{LangCF} and a $p$-adic interpretation of \eqref{hn2def}, 
they are able to prove an analogue of \eqref{iwasawa},
which  they call ``$(S,\{2\})$-Iwasawa theory'', 
and reads as follows (see \cite[Theorem 1.2 and \S4]{HK}). 

\begin{Theorem} 
\label{Theohn2}
There exist constants $m',\lambda'\ge0$ and $c'$ such that
\begin{equation}
\notag
\ordp\,h_{n,2}^{-}=m'p^n+\lambda' n+c'\qquad(n\gg0).
\end{equation}
\end{Theorem}
In the same article (see \cite[Remark 4.5]{HK}) the referee pointed out
that in Theorem \ref{Theohn2}
actually $m'=0$, this is, the growth of $\ordp\,h_{n,2}^{-}$ is linear, in analogy to  
the growth of $\ordp\,h_{n}^{-}$ given by Theorem \ref{Theohn}. The proof given by
the referee follows from the growth of the $p$-part
of some groups appearing in an exact sequence, which in turn comes from
an idelic interpretation of the $(S,T)$-ideal class group.

The aim of this article is to give a new proof of the fact that $m'=0$ in Theorem \ref{Theohn2}. 
The main idea is to relate the $p$-adic valuations of the numbers $h_{n,2}^{-}$ and $h_{n}^{-}$
to be able to use Theorem \ref{Theohn}. This allows us to obtain an explicit relation between the constants 
$\lambda'$, $c'$ above and the constants $\lambda$, $c$ appearing
in Theorem \ref{Theohn}. More precisely, we will prove the following.

\begin{Theorem} 
\label{TheoMain}
There exist constants $\lambda'\ge0$ and $c'$ such that
\begin{equation}
\notag
\ordp\,h_{n,2}^{-}=\lambda' n+c'\qquad(n\gg0).
\end{equation}
Moreover, these constants $\lambda'$ and $c'$ are related with the constants 
$\lambda$ and $c$ of Theorem \ref{Theohn} by means of $\lambda'=\lambda+\delta$
and $ c'=c+\delta$, where $\delta\ge0$ is an integer which does not depend on $n$ .
\end{Theorem}

The proof of Theorem \ref{TheoMain} will be given in \S\ref{proof}. 
In \S\ref{hn2hn} we relate the numbers $h_{n,2}^{-}$ and $h_{n}^{-}$ to be able to
use Theorem \ref{Theohn} and formula \eqref{hn2def}. In doing so, it appears a product
involving the values of some Dirichlet characters evaluated at $2$.
 In \S\ref{product} we compute the $p$-adic valuation
of this product (in a slightly more general form),
which we then use to compute $\ordp\,h_{n,2}^{-}$. 
Our proof is independent of Theorem \ref{Theohn2}.

Finally, we would like to mention that the study of more general functions 
than $L_E(s,\chi)$, and their $p$-adic analogues, was already done by Morita in \cite{MorL}.

\section{$h_{n,2}^{-}$ in terms of $h_{n}^{-}$}
\label{hn2hn}

First, for us to depend on the strength of Theorem \ref{Theohn}, we need to relate the 
values $L_E(0,\chi)$ and $L(0,\chi)=-B_{1,\chi}$. The methods we will use
are standard and we reproduce them here.

\begin{Lemma} 
\label{LemmaL}
We have the following identity:
\begin{equation}
\notag 
\prodchi L_E(0,\chi)=2^{\varphi(p^{n+1})/2}\prodchi\left(1-2\chi(2)\right)\prodchi B_{1,\chi}.
\end{equation}
\end{Lemma}
\begin{proof} Recall that $L_E(s,\chi)$ is defined by the series in \eqref{LEdef} for $\re(s)>0$.
Since we also need to work with the series of $L(s,\chi)$, we need the restriction $\re(s)>1$. Then
\begin{align}
\notag
2L(s,\chi)+L_E(s,\chi)&=2\sum_{k\ge1}\frac{\chi(k)}{k^s}+2\sum_{k\ge1}(-1)^k\frac{\chi(k)}{k^s}=
2\sum_{k\ge1}\left(1+(-1)^k\right)\frac{\chi(k)}{k^s}\\
\notag &=2\sum_{\substack{k\ge1\\ k\,\text{even}}}2\,\frac{\chi(k)}{k^s}=
4\sum_{j\ge1}\frac{\chi(2j)}{(2j)^s}=2\,\frac{\chi(2)}{2^{s-1}}\sum_{j\ge1}\frac{\chi(j)}{j^s}
=2\,\frac{\chi(2)}{2^{s-1}}L(s,\chi),
\end{align}
which gives $L_E(s,\chi)=2\left(\chi(2)/2^{s-1}-1\right)L(s,\chi)$.
By analytic continuation (see \cite[\S1]{MorL}), this is valid for $s=0$, and using \eqref{lzero}, we obtain
\begin{equation}
\notag
L_E(0,\chi)=2\left(2\chi(2)-1\right)(-B_{1,\chi})=2\left(1-2\chi(2)\right)B_{1,\chi}.
\end{equation}
Our identity then follows by taking the product 
over all odd primitive Dirichlet characters of conductor dividing $p^{n+1}$.
There are exactly $\varphi(p^{n+1})/2$ of them, which gives us the exponent in $2$.
\end{proof}

Therefore, we obtain $h_{n,2}^{-}$ in terms of $h_{n}^{-}$.

\begin{Proposition} 
\label{Relh}
The numbers $h_{n,2}^{-}$ and $h_{n}^{-}$ are related by means of
\begin{equation}
\notag
h_{n,2}^{-}=\frac{h_{n}^{-}}{p^{n+1}}\prodchi\left(1-2\chi(2)\right).
\end{equation}
\end{Proposition}
\begin{proof} This is a straightforward calculation using
\eqref{hnbn}, \eqref{hn2def} and Lemma \ref{LemmaL}.
\end{proof}

\section{Computation of $\,\ordp\prod_{\chi\,\text{odd}}\left(1-q\chi(q)\right)$}
\label{product}

Let $q$ be a prime number other than the fixed odd prime $p$. 
In this section we will compute the $p$-adic valuation of the product 
\begin{equation}
\label{tempq2}
\prodchi\left(1-q\chi(q)\right),
\end{equation}
which, in the case $q=2$, gives the product 
appearing in Proposition \ref{Relh}.

We shall write: $e=$ the multiplicative order of $q$ modulo $p$,
$f_n=$ the multiplicative order of $q$ modulo $p^{n+1}$,
and $m=\ordp\left(1-q^{e}\right)$. 

\begin{Lemma}
\label{lemmabasic}
Let the notation be as above, and suppose $n>m$. 
\begin{itemize}
\item[(i)] The multiplicative order of $q$ modulo $p^{n+1}$ is $f_n=ep^{n+1-m}$.
In particular,  $f_n$ and $e$ have
the same parity, and the parity of $f_n$ does not depend on $n$.
\item[(ii)]  We have that $\ordp\left(1-q^{f_n}\right)=n+1.$
\item[(iii)] If $f_n$ (or $e$) is even, then $\ordp\left(1+q^{f_n/2}\right)=n+1.$
\end{itemize}
\end{Lemma}
\begin{proof} We will use the following elementary result (see \cite[Theorem 3.6]{NathEMNT}):

\begin{center}
{\it Let $p$ be an odd prime and let $a\ne\pm1$ be an integer not divisible by $p$. Let $e$ 
be the multiplicative order of $a$ modulo $p$, and let 
$m=\ordp(1-a^e)$. Then, for $k\ge m$, the multiplicative order of $a$
modulo $p^{k}$ is $ep^{k-m}$.}
\end{center}
Also, we will use that $\ordp(1-a^{ep^{k-m}})=k$. To prove this, write
$a^e=1+\beta p^m$, where $\beta$ is not divisible by $p$. Then 
$$
\frac{a^{ep^{k-m}}-1}{p^k}=\frac{(1+\beta p^m)^{p^{k-m}}-1}{p^k}=
\beta+\beta' p^m,
$$
where $\beta'$ is an integer not divisible by $p$. Since $p$ does not divide $\beta$, this means
that $\ordp(1-a^{ep^{k-m}})=k$.

Now, (i) and (ii) follow immediately letting $a=q$, $k=n+1$ and $f_n=ep^{n+1-m}$ in
the above results. Recall that $p$ is odd,
hence $e$ and $f_n$ have the same parity.
In the case that $f_n$ (or $e$) is even, we can write
\begin{equation}
\notag
\ordp\left(1+q^{f_n/2}\right)=\ordp\left(\frac{1-q^{f_n}}{1-q^{f_n/2}}\right)=
\ordp\left(1-q^{f_n}\right)-\ordp\left(1-q^{f_n/2}\right).
\end{equation}
Since $q^{f_n/2}\equiv-1\,(\Mod p^{n+1})$ and since $p$ is odd, it follows that 
$\ordp\left(1-q^{f_n/2}\right)=0$. Hence, $\ordp\left(1+q^{f_n/2}\right)=
\ordp\left(1-q^{f_n}\right)=n+1$,
which proves (iii).
\end{proof}

Now we can compute the $p$-adic valuation of the product \eqref{tempq2}.

\begin{Proposition}
\label{propordq} 
Let the notation be as above, and suppose $n>m$. Then there exists an integer $d_q\ge1$, 
which is independent of $n$, such that
\begin{equation}
\notag
\ordp\prodchi\left(1-q\chi(q)\right)=(n+1)d_q.
\end{equation}
\end{Proposition}
\begin{proof} 
We will use the following result about characters of finite abelian groups: 

\begin{center}
{\it Let $A$ be a finite abelian group of order $N$ and let
$\hat A$ be its dual group. Let $a\in A$ have order $h$. Then 
$
\prod_{\chi\in\hat A}(1-\chi(a)T)=(1-T^h)^{N/h}.
$
}
\end{center}
We could not find an explicit reference for this result in the literature, so we give here a short proof,
courtesy of Keith Conrad:
Since $a$ has order $h$, the mapping $\hat A\to\mu_h$ given by 
$\chi\mapsto\chi(a)$ is a surjective homomorphism, so each $h$-th root of unity is a value
$N/h$ times. Thus the product is 
$\prod_{\chi\in\hat A}(1-\chi(a)T)=\prod_{z\in\mu_h}(1-zT)^{N/h}=(1-T^h)^{N/h}.$

Applying this for $A=(\Z/p^{n+1}\Z)^{\times}$ and $a=$ the class of $q$ in $A$, we obtain that
\begin{equation}
\label{aux1char}
\prod_{\chi\,\mathrm{mod}\,p^{n+1}}(1-\chi(q)T)=(1-T^{f_n})^{\varphi(p^{n+1})/f_n}.
\end{equation}
Also, applying the same result for $A=(\Z/p^{n+1}\Z)^{\times}/\{\pm1\}$ and $a=$ the class of $q$ in $A$, 
we obtain that
\begin{equation}
\label{aux2char}
\prodchieven(1-\chi(q)T)=
\begin{cases}
\left(1-T^{f_n/2}\right)^{\varphi(p^{n+1})/f_n}&\text{if $f_n$ (or $e$) is even},\\
\left(1-T^{f_n}\right)^{\varphi(p^{n+1})/(2f_n)}&\text{if $f_n$ (or $e$) is odd},
\end{cases}
\end{equation}
which is non-zero if $T=q$. Hence, letting $T=q$ and diving the product in \eqref{aux1char}
by the product in \eqref{aux2char},
we obtain that
\begin{equation}
\notag
\prodchi\left(1-q\chi(q)\right)=
\begin{cases}
\left(1+q^{f_n/2}\right)^{\varphi(p^{n+1})/f_n}&\text{if $f_n$ (or $e$) is even},\\
\left(1-q^{f_n}\right)^{\varphi(p^{n+1})/(2f_n)}&\text{if $f_n$ (or $e$) is odd}.
\end{cases}
\end{equation}
Now, define  $d_q=\varphi(p^{n+1})/f_n$ if $e$ is even, and $d_q=\varphi(p^{n+1})/(2f_n)$ if $e$ is odd.
By Lemma \ref{lemmabasic}.(i), we have that $\varphi(p^{n+1})/f_n=(p-1)p^{m-1}/e$.
Hence, $d_q$ does not depend on $n$, and we obtain that
\begin{equation}
\notag
\prodchi\left(1-q\chi(q)\right)=
\begin{cases}
\left(1+q^{f_n/2}\right)^{d_q}&\text{if $f_n$ (or $e$) is even},\\
\left(1-q^{f_n}\right)^{d_q}&\text{if $f_n$ (or $e$) is odd}.
\end{cases}
\end{equation}
Taking the $p$-adic valuation at both sides of this equation, our result follows from (ii) and (iii) in Lemma \ref{lemmabasic}.
\end{proof}

\section{Proof of Theorem \ref{TheoMain}}
\label{proof}

The proof now follows easily. Computing $\ordp$ at 
both sides of the identity given in Proposition \ref{Relh} we obtain that
\begin{equation}
\notag 
\ordp\,h_{n,2}^{-}=\ordp\,h_{n}^{-}-(n+1)+\ordp\prodchi\left(1-2\chi(2)\right).
\end{equation}
Using Proposition \ref{propordq} with $q=2$, and writing $\delta=d_2-1$, this becomes
\begin{equation}
\label{teo1}
\ordp\,h_{n,2}^{-}=\ordp\,h_{n}^{-}+(n+1)\delta\qquad(n>m),
\end{equation}
where the integer $\delta=d_2-1\ge0$ is given explicitly and does not depend on $n$. 
From Theorem \ref{Theohn},
there exist constants $\lambda\ge0$ and $c$ such that
\begin{equation}
\label{teo2}
\ordp\,h_n^{-}=\lambda n+c\qquad(n\gg0).
\end{equation}
Writing $\lambda'=\lambda+\delta\ge0$ and $c'=c+\delta$, 
equations \eqref{teo1} and \eqref{teo2} imply that
\begin{equation}
\notag
\ordp\,h_{n,2}^{-}=\lambda' n+c'
\end{equation}
for $n\gg0$ (such that also $n>m$), which proves Theorem \ref{TheoMain}.

\section*{Acknowledgments}

The author is very much in debt with the referee, who carefully read the preliminary version of this article 
and whose comments helped to improve the quality of the same. Also, the author would like to thank Keith Conrad
for providing the nice and short proof of the auxiliary result used in Proposition \ref{propordq}.
This research was supported by the author's CONICYT Chilean Doctoral Grant 2016.

\bibliographystyle{amsplain}
\bibliography{references}
\end{document}